\documentclass[12pt]{article}
\usepackage{amsfonts}
\usepackage{amssymb}
\usepackage{amsfonts,amssymb, amsmath, latexsym}
\usepackage{pst-all}
\usepackage{CJK}
\usepackage{mathrsfs}

\usepackage[left,pagewise,mathlines]{lineno}

\newtheorem{thm}{Theorem}

\newtheorem{lem}[thm]{Lemma}

\newtheorem{Def}{Definition}[section]

\newtheorem{obs}[thm]{Observation}

\newenvironment{pf}[1][Proof]{\noindent\textbf{#1.} }{\hfill\rule{1mm}{2mm}}

\makeatletter \@addtoreset{equation}{section} \makeatother

\def\g{\gamma}

\def\ma{\mathscr}

\begin{document}

\title{The $p$-Domination Number of Complete Multipartite Graphs\thanks{The work was
supported by NNSF of China (No.10711233) and the Fundamental
Research Fund of NPU (No. JC201150)}}

\author{
{You Lu$^a$ \quad Jun-Ming Xu$^b$}\\
{\small $^a$Department of Applied Mathematics,}\\
{\small             Northwestern Polytechnical University,}\\
{\small             Xi'an Shaanxi 710072, P. R. China}\\
{\small             Email: luyou@nwpu.edu.cn}\\ \\
{\small $^b$Department of Mathematics,}\\
{\small             University of Science and Technology of China,}   \\
{\small             Wentsun Wu Key Laboratory of CAS,}\\
{\small             Hefei, Anhui, 230026, P. R. China}\\
{\small             Email: xujm@ustc.edu.cn}  \\
}
\date{}
\maketitle

\begin{abstract}
Let $G=(V,E)$ be a graph and $p$ be a positive integer. A subset
$S\subseteq V$ is called a $p$-dominating set of $G$ if every vertex
not in $S$ has at least $p$ neighbors in $S$. The $p$-domination
number $\g_p(G)$ is the minimum cardinality of a $p$-dominating set
in $G$. In this paper, we establish an exact formula of the
$p$-domination number of all complete multipartite graphs for arbitrary
positive integer $p$.
\end{abstract}




\section{Induction}

For notation and graph-theoretical terminology not defined here, we refer the reader to \cite{bm08}.
Let $G=(V,E)$ be a finite
simple graph with vertex set $V=V(G)$ and edge set $E=E(G)$.
The {\it neighborhood} and {\it degree} of a vertex $v\in V$ are
$N_{G}(v)=\{u\in V : uv\in E\}$ and $d_G(v)=|N_G(v)|$, respectively.
 A {\it dominating set} of $G$ is a subset $S\subseteq V$ such that every
vertex of $V-S$ has at least one neighbor in $S$.
 The {\it domination number} $\gamma (G)$ is the minimum cardinality of
all dominating sets in $G$. The domination is a classical concept in
graph theory. The early literature on the domination with related
topics is, in detail, surveyed in the two outstanding books by Haynes,
Hedetniemi and Slater~\cite{hs981,hs982}.

Fink and Jacobson \cite{fj851,fj852}
generalized the concept of dominating
set.
Let $p$ be a positive integer.
A subset $D\subseteq V$ is a {\it $p$-dominating set} of $G$ if
$|N_G(v)\cap D|\geq p$ for each $v\in V-D$. The {\it
$p$-domination number} $\gamma_p(G)$ is the minimum cardinality of
all $p$-dominating sets in $G$. A $p$-dominating set $D$ with $|D|=\g_p(G)$
is called a $\g_p$-set of $G$ (for short, $\gamma_p(G)$-set). For $S,
T\subseteq V$, $S$ $p$-dominate $T$ in $G$ if
$|N_G(v)\cap S|\geq p$ for each $v\in T-S$. Clearly, the
$1$-dominating set is the well-known dominating set in a graph $G$, and so
 $\gamma_1(G)=\gamma(G)$. By the definition of $p$-dominating set, the following observation is obvious.

\begin{obs}\label{o1}
Every $p$-dominating set contains all the vertices with degree at most $p-1$.
\end{obs}

The determination of the $p$-domination number for graphs seems to
be a difficult problem. In 1989, Jacobson and Peters \cite{jp89}
showed that the problem is NP-complete in general graphs. In 1994,
Bean, Henning and Swart \cite{bhs94} proved the problem remains
NP-complete in bipartite or chordal graphs. These results show that
the following study is of important significance.
\begin{itemize}
\item Find the lower and upper bounds of $\g_p$ with difference as small as possible.
\item Determine exact values of $\g_p$ for some graphs, specially well-known networks.
\end{itemize}

Many works focused on the bounds of $\g_p$ for general graphs or
some special classes of graphs (see, for example,
\cite{bcv06,cr90,dghpv,f85,rv07}). Very recently, Chellali et al.
\cite{cfhv11} have given an excellent survey on this topics. Until
now, however, no research has been done on calculating the exact
values of $\g_p$ even for some particular graphs except \cite{s09}.
In \cite{s09}, the author obtained the exact $2$-domination number
of the toroidal grid graphs $C_m\square C_n$ in some cases.

In this paper, we give an exact formula of $\g_p$ for arbitrary
positive integer $p$ and the complete $t$-partite graph
$K_{n_1,n_2,\cdots,n_t}$.

Through this paper, the graph $G$ always denotes a complete
$t$-partite graph $K_{n_1, n_2, \cdots, n_t}$ with $t$-partition
$\{V_1,V_2,\cdots, V_t\}$, $N_t=\{1,2,\cdots,t\}$ and
 $$
 f(I)=\sum_{i\in I}n_i\ \ {\rm for}\ I\subseteq N_t.
 $$
Note that if $t=1$ or $t=2$ and $f(N_t)\leq p$ then $\g_p(G)=|V(G)|$
by Observation \ref{o1}. Thus, we always assume $t\geq 2$ and
$f(N_t)> p$.

\section{Optimal $\g_p$-sets of $G$}
For any $D\subseteq V(G)$, define
 \begin{eqnarray*}
 D_i=V_i\cap D\ \ {\rm for\ each}\ i\in N_t, \mbox{\ \ and\ \ }
 I_D = \{i\in N_t:\  |D_i|=|V_i|\}.
 \end{eqnarray*}

\begin{lem}\label{lem2}
If $t\geq 2$ and $f(N_t)> p$, then
 $$
 \g_p(G)\leq \min\{f(I) : I\subseteq N_t \mbox{ with $f(I)\geq p$}\}.
 $$
with equality if $G$ has a $\g_p$-set $D$ with $f(I_D)\geq p$.
\end{lem}

\begin{pf}
Let $I\subseteq N_t$ with $f(I)\geq p$ and $S=\bigcup_{i\in I}V_i$.
Then
  $$
  |S|=\sum_{i\in I}|V_i|=\sum_{i\in I}n_i=f(I)\geq p.
  $$
Since $G$ is a complete $t$-partite graph, for any $v\in V(G)-S$, we
have $S\subseteq N_G(v)$ and so $|N_G(v)\cap S|=|S|\geq p$. This
implies that $S$ is a $p$-dominating set of $G$, and so
 $$
 \g_p(G)\leq \min\{f(I) : I\subseteq N_t \mbox{ with $f(I)\geq p$}\}.
 $$

On the other hand, let $D$ be a $\g_p(G)$-set with $f(I_D)\geq
p$. Then
 \begin{eqnarray*}
  \g_p(G)=|D|&\geq& \sum\limits_{i\in I_D}|V_i|=\sum\limits_{i\in I_D}n_i=f(I_D)\\
  &\geq& \min\{f(I) : I\subseteq N_t \mbox{ with $f(I)\geq p$}\}.
   \end{eqnarray*}

 The lemma follows.
\end{pf}

\begin{lem}\label{lem3}
If $t\geq 2$ and $f(N_t)> p$, then $|I_D|\leq t-2$ for any
$\g_p(G)$-set $D$ with $f(I_D)< p$.
\end{lem}

\begin{pf}
Clearly $|I_D|\le t-1$ by $f(N_t)> p>f(I_D)$. If $|I_D|=t-1$, then
there is a unique index $i_0\in N_t$ such that $N_t-I_D=\{i_0\}$. By
the definition of $I_D$, $V(G)-V_{i_0}\subseteq D$ and there exists
a vertex $x$ in $V_{i_0}$ but not in $D$.  Since $D$ is a $\g_p(G)$-set and $f(I_D)< p$, we can deduce a contradiction as follows:
 $$
 p\leq |N_G(x)\cap D|=|V(G)-V_{i_0}|=\sum_{i\in I_D}n_i=f(I_D)< p.
 $$
Hence $|I_D|\leq t-2$.
\end{pf}

\vskip6pt


For a $\g_p(G)$-set $D$ with $|I_D|<t$, $|D|=f(I_D)$
$+\sum\limits_{i\in N_t-I_D}|D_i|$. By Lemma~\ref{lem3}, the value
of $\left|\,|D_i|-\frac{|D|-f(I_D)}{t-|I_D|}\right|$ is well-defined
for any $i\in N_t-I_D$ if $t\geq 2$ and $f(N_t)> p$. Define
 $$
\mu(D)=\sum_{i\in
N-I_D}\left||D_i|-\frac{|D|-f(I_D)}{t-|I_D|}\right|.
 $$


\begin{Def}
A $\g_p(G)$-set $D$ is called to be {\it optimal} if the following
conditions hold: 1) $f(I_D) <p$; 2) $|I_D|\geq |I_S|$ for any
$\g_p(G)$-set $S$; 3) $\mu(D)\le \mu(S)$ for any $\g_p(G)$-set $S$
with  $I_D=I_S$.
\end{Def}

By the definition, if each $\g_p(G)$-set $D$ has $f(I_D)<p$, then
there must be at least one optimal $\g_p$-set in $G$. To obtain the
upper bound of $\g_p(G)$,
by Lemma \ref{lem2}, we only need to consider the case that
every $\g_p(G)$-set $D$ satisfies $f(I_D)<p$.
We investigate properties of optimal $\g_p$-sets starting with the
following critical lemma.

\begin{lem}\label{lem4}
$\left||D_i|-|D_j|\right|\leq 1$ for any optimal $\g_p(G)$-set $D$
and $i,j\in N_t-I_D$.
\end{lem}

\begin{pf}
 By Lemma \ref{lem3}, $t-|I_D|\geq 2$ and so $N_t-I_D\neq \emptyset$. Let
 $$
 |D_{s}|=\max\{|D_i|:\ i\in N_t-I_D\}\ \ {\rm and}\ \
 |D_{w}|=\min\{|D_i|:\ i\in N_t-I_D\}.
 $$

Suppose, to be contrary, that $|D_{s}|-|D_{w}|\geq 2$. Clearly,
$|D_{s}|\geq 2$. Since $w\in N_t-I_D$, $D_{w}\varsubsetneq V_{w}$. Hence there are
$x\in D_{s}$ and $y\in V_{w}-D_{w}$. Let
  $$
  D^*=(D-\{x\})\cup \{y\}.
  $$
Then
\begin{eqnarray}\label{e2.4}
     I_{D^*}=\left\{\begin{array}{ll}
              I_D & \mbox{ if $|D_w|<|V_w|-1$};\\
              I_D\cup \{w\} & \mbox{ if $|D_w|=|V_w|-1$}.
              \end{array}
              \right.
\end{eqnarray}
Thus $I_D\subseteq I_{D^*}$. We first claim that $D^*$ is a
$\g_p(G)$-set. In fact, it is easy to see that $D^*$ can
$p$-dominate $V(G)-V_{w}$. By the choice of $s\in N_t-I_D$,
$V_{s}-D_s\neq \emptyset$. Since $D-D_s$ can $p$-dominate
$V_{s}-D_s$, we have $|D|-|D_{s}|\geq p$. It follows that, for any
vertex $z\in V_{w}-D^*$,
 \begin{eqnarray*}
|N_G(z)\cap D^*|=|D^*|-|D^*_w|=|D|-(|D_{w}|+1)\geq |D|-|D_{s}|+1\geq
p+1,
\end{eqnarray*}
which means that $D^*$ can $p$-dominate $z$ and, hence, $D^*$ is a
$\g_p(G)$-set.

By the second condition of the optimality of $D$, we have $|I_D|\geq |I_{D^*}|$. Thus $I_D=I_D^*$ by $I_D\subseteq I_{D^*}$. Combined with $|D|=|D^*|=\g_p(G)$, we can obtain that
$$\frac{|D|-f(I_D)}{t-|I_D|}
 =\frac{|D^*|-f(I_{D^*})}{t-|I_{D^*}|}.$$
For convenience, we use the notation $\lambda$ to represent them.

We now show $\mu(D^*)-\mu(D)<0$.
Since $|D|=\g_p(G)=f(I_D)+\sum_{i\in N_t-I_D}|D_i|$,
 \begin{eqnarray*}
 \lambda=\frac{1}{t-|I_D|}\sum_{i\in N_t-I_D}|D_i|.
 \end{eqnarray*}
By $|D_{s}|-|D_{w}|\geq 2$ and the choices of $s$ and $w$, we have
that
 $$
 |D_{w}|+1\leq |D_{s}|-1 \mbox{\ \ and\ \ } |D_{w}|< \lambda< |D_{s}|.
 $$
It follows that
\begin{eqnarray*}
\mu(D^*)-\mu(D)&=&\sum_{i\in N_t-I_{D^*}}|(|D^*_i|-\lambda)|-\sum_{i\in N_t-I_D}|(|D_i|-\lambda)|\\
              &=&|(|D^*_w|-\lambda)|+|(|D^*_s|-\lambda)|-(\lambda-|D_{w}|)-(|D_{s}|-\lambda)\\
              &=&|(|D_{w}|+1-\lambda)|+|(|D_{s}|-1-\lambda)|-(|D_{s}|-|D_{w}|)\\
              &=&\left\{\begin{array}{ll}
                             2(|D_{w}|-\lambda)& \mbox{\ \ if\ \ $\lambda< |D_{w}|+1$}\\
                             -2                  &  \mbox{\ \ if\ \ $|D_{w}|+1\leq \lambda\leq |D_{s}|-1$}\\
                             2(\lambda-|D_{s}|)&   \mbox{\ \ if\ \ $\lambda> |D_{s}|-1$}\\
                         \end{array}
              \right.\\
              &<& 0.
\end{eqnarray*}
This contradicts with the third condition of the optimality of $D$, and so
$|D_{s}|-|D_{w}|\leq 1$.

The lemma follows.
\end{pf}

\vskip6pt

For an optimal $\g_p(G)$-set $D$, $t-|I_D|\geq 2$ by Lemma \ref{lem3},  and so $N_t-I_D\neq \emptyset$.  Thus we denote
 \begin{equation}\label{e2.1}
 k=\max\{|D_i|:\ i\in N_t-I_D\}\ \ {\rm and}\  \
 \ell=\min\{|D_i|:\ i\in N_t-I_D\}.
 \end{equation}
If $k\ne \ell$, then $k=\ell+1$ by Lemma~\ref{lem4}. Define
 \begin{equation}\label{e2.2}
 \begin{array}{rl}
  & A=\left\{\begin{array}{ll}
 \{i\in N_t-I_D: |D_i|=\ell+1\}\ &\ {\rm if}\ k=\ell+1;\\
 \emptyset \ &\ {\rm if}\ k=\ell,
  \end{array}\right.\\
  &B=\{i\in N_t-I_D: |D_i|=\ell\}.
   \end{array}
 \end{equation}
Then $\{A,B\}$ is a partition of $N_t-I_D$ and $B\neq \emptyset$.

\begin{lem}\label{lem5}
$|A|=0$ or $2\leq |A|\leq t-|I_D|-1$ for any optimal $\g_p(G)$-set
$D$.
\end{lem}

\begin{pf}
Since $\{A,B\}$ is a partition of $N_t-I_D$ and $B\neq \emptyset$,  it is obvious that $|A|\leq t-|I_D|-1$.
We now show $|A|\ne 1$. Assume to the contrary that $|A|=1$.

Let $A=\{i_1\}$. Then $|D_{i_1}|=\ell+1\geq 1$ and $V_{i_1}-D_{i_1}\neq \emptyset$ since $i_1\in
N_t-I_D$. Since $D-D_{i_1}$ $p$-dominates $V_{i_1}-D_{i_1}$, we have
$|D|-|D_{i_1}|\geq p$. Take any vertex $x\in D_{i_1}$ and let
$$D'=D-\{x\}.$$
 Consider any vertex $y$ in $V(G)-D'$. If $y\in V_{i_1}$, then
$$|N_G(y)\cap D'|=|D'|-|D'_{i_1}|=(|D|-1)-(|D_{i_1}|-1)|=|D|-|D_{i_1}|\geq p.$$
If $y\notin V_{i_1}$, then there exists some $j\in B$ such that
$y\in V_j$. Noting $|D_j|=|D_{i_1}|-1$, we have that
$$
|N_G(y)\cap D'|=|D'|-|D'_j|=(|D|-1)-|D_j|=|D|-|D_{i_1}|\geq p.
$$
Hence $D'$ is a $p$-dominating set of $G$ with
$|D'|=|D|-1=\g_p(G)-1$, a contradiction. The lemma follows.
\end{pf}

\begin{lem}\label{lem6}
$\g_p(G)\geq p+\ell+\delta_A$ for any optimal $\g_p(G)$-set $D$,
where $\ell$ and $A$ are defined in (\ref{e2.1}) and (\ref{e2.2}),
respectively, $\delta_A$ is the characteristic function on $A$,
i.e., $\delta_A=0$ if $|A|=0$ and $\delta_A=1$ otherwise.
\end{lem}

\begin{pf}
Note that $N_t-I_D\ne\emptyset$ and $V_i-D_i\neq \emptyset$ for $i\in
N_t-I_D$. To $p$-dominate $V_i-D_i$,
 $
|D-D_i|= |D|-|D_i|\geq p
 $
for $i\in N_t-I_D$.

If $|A|=0$, then $\delta_A=0$ and $N_t-I_D=B$. For any $i\in B$,
$|D_i|=\ell$ by (\ref{e2.2}), and so $\g_p(G)=|D|\geq
p+|D_i|=p+\ell=p+\ell+\delta_A$.

If $|A|\neq 0$, then $\delta_A=1$.
For $i\in A$, $|D_i|=\ell+1$ by (\ref{e2.2}). Thus
$\g_p(G)=|D|\geq p+|D_i|=p+\ell+1=p+\ell+\delta_A$.

The lemma follows.
\end{pf}

\begin{lem}\label{lem7}
 $\lceil\frac{p-f(I_D)}{t-|I_D|-1}\rceil\leq n_i$
 for any optimal $\g_p(G)$-set $D$ and $i\in N_t-I_D$.
\end{lem}

\begin{pf}
Let $N-I_D=A\cup B$ as defined in (\ref{e2.2}). Then $|D_i|=\ell+1$ for $i\in A$ and $|D_j|=\ell$ for
$j\in B$. Note that $n_i=|V_i|\geq |D_i|+1\geq \ell+1$ for any $i\in N_t-I_D$. It follows that
 \begin{eqnarray*}
 |D|&=&f(I_D)+f(A)+f(B)\\
           &=&f(I_D)+|A|(\ell+1)+(t-|I_D|-|A|)\ell\\
           &=&f(I_D)+(t-|I_D|-1)\ell+\ell+|A|,
           \end{eqnarray*}
from which we have
 $$
 \begin{array}{rl}
\lceil\frac{p-f(I_D)}{t-|I_D|-1}\rceil&=\ell+\lceil\frac{|A|-\delta_A}{t-|I_D|-1}-\frac{|D|-(p+\ell+\delta_A)}{t-|I_D|-1}\rceil\\
                                      &\leq\ell+\lceil\frac{|A|-\delta_A}{t-|I_D|-1}\rceil \hspace{1cm} \mbox{(by Lemma~\ref{lem6})}\\
                                      &\leq \ell+\delta_A\hspace{2.3cm} \mbox{(by Lemma~\ref{lem5})}\\
                                      &\leq \ell+1 \\
                                      &\leq n_i\ \ \ {\rm for\ any}\ i\in N_t-I_D\\
                                      \end{array}
                                      $$
as desired, and so the lemma follows.
\end{pf}

\section{Main results}

In this section, we will give an exact formula of $\g_p$ for a
complete $t$-partite graph $G=K_{n_1,n_2,\cdots,n_t}$. By
Lemma~\ref{lem2}, if $G$ contains a $\g_p$-set $D$ with
$f(I_D)\geq p$, then
 $$
 \g_p(G)=\min\{f(I) : I\subseteq N_t \mbox{ with $f(I)\geq p$}\}.
 $$
Thus, we only need to consider the case of $f(I_D)< p$ for any $\g_p(G)$-set $D$. In this case, $G$ must have optimal $\g_p(G)$-sets. Moreover, for any optimal
$\g_p(G)$-set $D$,  $|I_D|\leq t-2$ by Lemma~\ref{lem3}, and
$\lceil\frac{p-f(I_D)}{t-|I_D|-1}\rceil\leq n_i$ for any $i\in N_t-I_D$ by
Lemma~\ref{lem7}. Thus, the following family $\ma{I}_p$ of the
subsets of $N_t$ is well-defined.
 $$
 \begin{array}{rl}
\ma{I}_p=\{I\subset N_t: |I|\leq t-2, f(I)< p \mbox{ and }
                         \lceil\frac{p-f(I)}{t-|I|-1}\rceil\leq n_i \mbox{ for each } i\in N_t-I\}.
 \end{array}
 $$
Some examples of $\ma{I}_p$ for $G=K_{2,2,10,17}$ can be found in
Table \ref{Table1}.
\begin{table}[h]
  \begin{center}
    \begin{tabular}{|c|c|c|c|c|}\hline
      $p$      &  $s_1$   &     $\ma{I}_p$                      &       $s_2$   &      $\g_p(G)$         \\ \hline
      $1$      &  $2$     &  $\{\emptyset\}$                    &       $1$     &      $s_1=p+s_2=2$     \\ \hline
      $2$      &  $2$     &  $\{\emptyset\}$                    &       $1$     &      $s_1=2$           \\ \hline
      $3$      &  $4$     &  $\{\emptyset,\{1\},\{2\}\}$        &       $1$     &      $s_1=p+s_2=4$     \\ \hline
      $4$      &  $4$     &  $\{\emptyset,\{1\},\{2\}\}$        &       $1$     &      $s_1=4$           \\ \hline
      $5$      &  $10$    &  $\{\emptyset,\{1\},\{2\},\{1,2\}\}$&       $1$     &      $p+s_2=6$         \\ \hline
      $6$      &  $10$    &  $\{\emptyset,\{1\},\{2\},\{1,2\}\}$&       $2$     &      $p+s_2=8$         \\ \hline
      $7$      &  $10$    &  $\{\{1,2\}\}$                      &       $3$     &      $s_1=p+s_2=10$    \\ \hline
      $9$      &  $10$    &  $\{\{1,2\}\}$                      &       $5$     &      $s_1=10$          \\ \hline
      $11$     &  $12$    &  $\{\{1\},\{2\},\{3\},\{1,2\}\}$    &       $1$     &      $s_1=p+s_2=12$    \\ \hline
      $13$     &  $14$    &  $\{\{3\},\{1,2\},\{1,3\},\{2,3\}\}$&       $1$     &      $s_1=p+s_2=14$    \\ \hline
      $14$     &  $17$    &  $\{\{3\},\{1,2\},\{1,3\},\{2,3\}\}$&       $2$     &      $p+s_2=16$        \\ \hline
      $15$     &  $17$    &  $\emptyset$                        &       $\infty$&      $s_1=17$          \\ \hline
    \end{tabular}
  \end{center}
 \caption{Examples of $s_1$, $\ma{I}_p$, $s_2$ and $\g_p(G)$ for $G=K_{2,2,10,17}$, where $N_4=\{1,2,3,4\}$.}\label{Table1}
\end{table}


Let
$$
\begin{array}{rl}
  &s_1=\min\{f(I) : I\subseteq N_t \mbox{ with $f(I)\geq p$}\}\ \ {\rm
 and}\\
 &\\
 &s_2=\left\{  \begin{array}{ll}
                        \min\{\lceil\frac{p-f(I)}{t-|I|-1}\rceil : I\in \mathscr{I}_p\} & \mbox{if }\mathscr{I}_p\neq \emptyset;\\
                        \infty & \mbox{if }\mathscr{I}_p=\emptyset.
                      \end{array}
                  \right.
\end{array}
$$

\begin{lem}\label{lem8}
Let $G=K_{n_1, n_2, \cdots, n_t}$ with $t\geq 2$ and $f(N_t)> p$. Then
$
\g_p(G)\le p+s_2.
$
\end{lem}

\begin{pf}
If $\mathscr{I}_p=\emptyset$, then $s_2=\infty$ and so
$\g_p(G)<p+s_2$. Assume that $\mathscr{I}_p\neq \emptyset$ below.
Let $I\in \ma{I}_p$ (without loss of generality, say
$I=\{1,\cdots,k\}$) with
 $$
 \begin{array}{rl}
 k\leq t-2,\  f(I)< p
 \mbox{ and }s_2=\lceil\frac{p-f(I)}{t-k-1}\rceil\leq n_i\
 \mbox{for each } i\in \{k+1,\cdots,t\}.
 \end{array}
$$

Since $t-k-1\geq 1$ and $p-f(I)>0$, there are two integers $q$ and
$r$ with $q\geq 0$ and $0\leq r\leq t-k-2$ such that
$$
p-f(I)=q(t-k-1)+r.
$$
Then for each $i\in \{k+1,\cdots,t\}$,
 \begin{equation}\label{e1}
 n_i\geq s_2=\left\{\begin{array}{ll}
                q+1 & \mbox{ if $r\neq 0$};\\
                  q & \mbox{ if $r=0$}.
                 \end{array}
      \right.
 \end{equation}
Thus, we can choose $D\subseteq V(G)$ such that
$$D=(V_1\cup \cdots \cup V_k) \cup (V_{k+1}'\cup \cdots \cup V_{k+r}') \cup (V_{k+r+1}'\cup \cdots \cup V_{t-1}')\cup V_t',$$
where, for each $i\in \{k+1,\cdots,t\}$, $V'_i$ is a subset of $V_i$
satisfying
\begin{eqnarray} \label{e2}
       |V_i'|=\left\{\begin{array}{ll}
                        q+1 & \mbox{ if } k+1\leq i\leq k+r\\
                        q   & \mbox{ if } k+r+1\leq i\leq t-1\\
                        s_2    & \mbox{ if } i=t.
                    \end{array}
              \right.
\end{eqnarray}
Thus,
\begin{eqnarray*}
|D|&=&\sum_{i=1}^k|V_i|+\sum_{i=k+1}^{k+r}|V_i'|+\sum_{j=k+r+1}^{t-1}|V_j'|+|V_t'|\\
    &=&(n_1+\cdots+n_k)+r(q+1)+(t-k-r-1)q+s_2 \\
    &=&(f(I)+q(t-k-1)+r)+s_2\\
    &=&p+s_2.
\end{eqnarray*}

To complete the proof, we only need to show that $D$ is a
$p$-dominating set of $G$.  To this aim, let $v$ be any vertex in
$V(G)-D$. By the choice of $D$, there is some $i_0\in
\{k+1,\cdots,t\}$ such that $v\in V_{i_0}-V'_{i_0}$. Since $G$ is a
complete $t$-partite graph,
$$
|N_G(v)\cap D|=|D|-|V'_{i_0}|=p+s_2-|V'_{i_0}|.
$$
By (\ref{e1}) and (\ref{e2}), we have
\begin{eqnarray*}
s_2-|V'_{i_0}|
      &=& \left\{\begin{array}{ll}
                 1 & \mbox{ if $r\neq 0$ and $k+r+1\leq i_0\leq t-1$}\\
                 0 & \mbox{ otherwise}\\
                 \end{array}
      \right.\\
      &\geq& 0.
      \end{eqnarray*}
It follows that $|N_G(v)\cap D|=p+s_2-|V'_{i_0}|\geq p$, which
implies that $D$ can $p$-dominate $x$. Hence $D$ is a $p$-dominating
set of $G$.
The lemma follows.
\end{pf}

\vskip6pt

We now state our main result as follows.

\begin{thm}\label{thm9}
For any integer $p\geq 1$ and a complete $t$-partite graph
$G=K_{n_1, n_2, \cdots, n_t}$ with $t\geq 2$ and $f(N_t)> p$,
$$
\g_p(G)=\min\{s_1,p+s_2\}.
$$
\end{thm}

\begin{pf}
 From Lemmas \ref{lem2} and \ref{lem8}, we can obtain that $\g_p(G)\leq \min\{s_1,p+s_2\}$, and if $G$ has a $\g_p$-set $D$ with $f(I_D)\geq p$ then $\g_p(G)= s_1\geq \min\{s_1,p+s_2\}$.

In the following, assume that every $\g_p(G)$-set $D$ satisfies $f(I_D)< p$.
Let $D$ be an optimal $\g_p(G)$-set. To the end, we only need to show $\g_p(G)\geq p+s_2$.

Since $|I_D|\leq t-2$ by Lemma~\ref{lem3} and
$\lceil\frac{p-f(I_D)}{t-|I_D|-1}\rceil\leq n_i$ for any $i\in N_t-I_D$ by
Lemma~\ref{lem7}, we have $I_D\in \ma{I}_p$, and so
 $
 \begin{array}{rl}
\lceil\frac{p-f(I_D)}{t-|I_D|-1}\rceil\geq s_2.
\end{array}
$
From the proof of Lemma~\ref{lem7}, we know that $\ell+\delta_A\geq \lceil\frac{p-f(I_D)}{t-|I_D|-1}\rceil$. Hence, by Lemma \ref{lem6},
 $$
 \begin{array}{rl}
 \g_p(G)\geq p+\ell+\delta_A\geq p+\lceil\frac{p-f(I_D)}{t-|I_D|-1}\rceil\geq p+s_2.
 \end{array}
 $$
The theorem follows.
\end{pf}

Some illustrations of $s_1$, $s_2$ and $\g_p(G)=\min\{s_1,p+s_2\}$ for the complete $4$-partite graph $G=K_{2,2,10,17}$ are shown in Table \ref{Table1}.

\end{document}